\documentclass[12pt,a4paper]{article} %[blank_reccur-sgroup.tex 10.08.20]
\usepackage{amsfonts,amssymb} 

\usepackage[margin=25mm, a4paper]{geometry}  % journal quality

\usepackage{tikz}

\def\n{\noindent}  \def\?#1{}
\def\IZ{{\mathbb{Z}}}  
   \def\cB{{\cal B}} 
\def\cT{{\cal T}}   \def\cF{{\cal F}}   \def\cS{{\cal S}}  \def\bp{{\mathbf{p}}}
\def\bT{{\mathbf{T}}}   \def\Rec{{\rm Rec}}

 \def\supp{{\rm supp}}
  \def\t{\tilde} \def\bdelete#1{}
\def\mod1{\,({\rm mod\ } 1)\,}

\def\beq#1#2{\begin{equation} \label{#1} #2 \end{equation}}
 
\def\function#1{\left\{\!\!\!\begin{array}{ll} #1 \end{array} \right.}
\def\proof{\smallskip \noindent {\bf Proof. \ }}       %start of proof
\def\blanksquare{\,\,\,$\sqcup\!\!\!\!\sqcap$}         %blank  square
\def\qed{\hfill\blanksquare\linebreak\smallskip\par}   %end of proof

\def\thname{Theorem}  \def\lmname{Lemma}    \def\prname{Proposition}
\def\dfname{Definition}  \def\crname{Corollary}  \def\rmname{Remark}
\def\exname{Example}  

\newtheorem{theorem}{\thname}%[section]   %Numbering: Theorem--Other section
\newtheorem{lemma}{\lmname}%[section]     %{lemma}[theorem]{Lemma}   subsection
\newtheorem{corollary}[lemma]{\crname}   %lemma
\newtheorem{example}{\exname}%[section]

\newtheorem{dftn}{\dfname}%[section]
\newenvironment{definition}{\begin{dftn}\rm}{\end{dftn}} %section
\def\bdef#1{\begin{definition} #1 \end{definition}}
\newtheorem{rmrk}[lemma]{\rmname}
\newenvironment{remark}{\begin{rmrk}\rm}{\end{rmrk}}     %lemma

\makeatletter\def\fps@figure{htbp}\makeatother %figure pos: tbp - standard

\begin{document}

\title{%(DCDS-A) 
        Recurrence for measurable semigroup actions}
\author{Michael Blank\thanks{
        Institute for Information Transmission Problems RAS
        (Kharkevich Institute);}
        \thanks{National Research University ``Higher School of Economics'';
        e-mail: blank@iitp.ru}
       }
\date{August 10, 2020} %\today}
\maketitle

\begin{abstract}%
We study qualitative properties of the set of recurrent points of finitely 
generated free semigroups of measurable maps. In the case of a single 
generator the classical Poincare recurrence theorem shows that these 
properties are closely related to the presence of an invariant measure. 
Curious, but otherwise it turns out to be possible that almost all points are 
recurrent, while there is an wandering set of positive (non-invariant) measure.
For a general semigroup the assumption about the common invariant 
measure for all generators looks somewhat unnatural (despite being 
widely used). Instead we give abstract conditions (of conservativity type) 
for this problem and propose a weaker version of the recurrent property. 
Technically,  the problem is reduced to the analysis of the recurrence 
of a specially constructed Markov process. Questions of inheritance of 
the recurrence property from the semigroup generators to the entire 
semigroup and vice versa are studied in detail and we demonstrate 
that this inheritance might be rather unexpected.
\end{abstract}%

{\small\n
}

%%%%%%%%%%%%%%%%%%%%%%%%%%%%%%%%
\section{Introduction}

Recall that a (discrete time) semigroup of maps $\cT$ is a closed under 
composition collection of maps acting on the same topological space $X$. 
One says that $\cT$ is finitely generated if there is a finite number 
of elements $\bT:=\{T_1,T_2,\dots,T_d\}\subset\cT$ whose compositions 
generate $\cT$. Given a point $x\in X$ one can define its trajectory under 
the action of the semigroup $\cT$ in two very different ways. First of them, 
to which we will refer as a global trajectory, is the union of all images of 
the point $x$ under the action of the maps from $\cT$, i.e. $\cup_{T\in\cT}Tx$. 
The main disadvantage of this definitions is that the points of the global 
trajectory are not ordered. To avoid this difficulty we consider another approach, 
to which we will refer as a realization (as in random processes), and which 
can be described as an infinite sequence of compositions of the generators 
of the semigroup evaluated at the point $x$, i.e. 
$\dots\circ T_{i_n}\circ T_{i_{n-1}}\circ\dots\circ T_{i_1}x$. The union of all 
realizations is exactly the global trajectory of the point under the action of the 
semigroup. 

The idea of recurrence is very old and simple. Roughly speaking it means that 
a trajectory (or a realization of a random process) returns eventually arbitrary 
close to its initial point. From the point of view of dynamical systems theory, 
recurrence is a distant generalization of the notion of periodicity. 
In terms of the global trajectories this means that there is a sequence 
$\{\tau_{j_1},\dots,\tau_{j_n},\dots\}$ of elements $\tau_j\in\cT$ such that 
$\lim_{n\to\infty}\tau_{j_n}x=x$, while the more natural ``realization'' approach 
indeed describes the returning of the consecutive images of the point $x$ along 
the realization to any neighborhood of $x$. 

In distinction to the most known scenario -- random walk on a discrete phase space 
(see e.g. \cite{Fe,MT,Or,Sp}), in the general case one cannot expect the returning exactly 
to the initial position. We introduce several types of the recurrence property 
(see bellow) and study conditions under which there are recurrent points of each 
of the types and when the set of such points is large enough (say of positive reference 
measure, which needs not to be dynamically invariant). 

The simplest (and the most known) instance of the finitely generated semigroup 
is a discrete time dynamical system, which corresponds to the case of a single 
generator. In this case the points of a realization $\{T^nx\}_{n\in\IZ_+}$ are naturally 
ordered with respect to time $n$. In the case when the map $T$ is continuous 
the situation with the recurrence is relatively simple and a number of results 
are already known (see \cite{Ni,MF}). On the other hand, in the measurable category 
the situation is more complicated and the constructions described in Section~\ref{s:S1} 
are new even in the case of a single generator. 

Technically all our results related to the abundance of recurrent points are based in 
one way or another on the celebrated Poincare Theorem (see e.g. \cite{Si,Fu1}), 
which asserts that for a $\mu$-measure preserving dynamical system $\mu$-a.a. 
points return under dynamics to any given measurable set. If the map is continuous 
this result immediately implies the recurrence of $\mu$-a.a. points in the phase 
space. However for a measurable category of maps that we consider here 
some additional care is necessary to prove this claim (see \cite{Bl2019}). 
Curiously, it turns out to be possible that almost all points are recurrent, 
while there is an wandering set of positive (non-invariant) measure.

Observe that the Poincare Theorem gives information about returns to a given set 
in terms of a dynamically invariant measure\footnote{Up to pure technical generalizations 
   like incompressibility ($A\in\cB$ with $A\subseteq T^{-1}A$ implies $m( T^{-1}A\setminus A)=0$), 
   see, e.g. \cite{Ha, Kr}. The latter property is non-constructive 
   and practically unverifiable, unlike our condition (\ref{e:poincare}) which allows 
   a simple check. }, 
whose support may be quite small and which in general might not exist. 
In order to deal with the latter case and to study recurrent properties of points 
not belonging to the support of the invariant measure we are developing conditions 
(of conservativity type) under which a (not dynamically invariant) reference 
measure (say Lebesgue measure) plays the same role. For Markov chains 
this was partially done in our previous work \cite{Bl2019} and here we extend this 
approach for the case of finitely generated semigroups of measurable maps. 
It is worth noting that the situation with the recurrence property for 
semigroups is much more involved compared to the case of a single self-map 
and one needs to develop new technical tools to its analysis.

We already noted that the absence of a common invariant measure (CIM) is a serious 
obstacle in the study of ergodic properties of a semigroup. The absence of an 
invariant measure being exotic for a single self-map in the compact phase space setting  
(see Section~\ref{s:S1}) is generic even in the case of a semigroup with 2 generators.  
Under a rather restricting and non-generic condition of the presence of CIM and assuming 
that the semigroup consists only of continuous maps a number of results are known  
about qualitative recurrence property, e.g. of Kac's Lemma type claims (see \cite{CRV,Fu2,Yo} 
and further references therein). However their counterparts for the measurable category 
of self-maps, not speaking about the absence of CIM, are not known at present. 
Let us mention also a version of the Khinchine type recurrence theorem studied 
in the case of a single map in \cite{Bo,BFW} and generalized for continuous 
group actions with CIM in \cite{Ki}. How to deal with this problem without specifying 
CIM is one of interesting questions for future research.

The paper is organized as follows. In Section~\ref{s:S1} we study the case of a 
semigroup with a single generator -- a measurable dynamical system. It turns out that 
even in this simplest setting there are quite a few properties of the recurrent points 
not considered earlier. In particular, we introduce and study a weak\footnote{As well as 
    two stronger versions.} version of the recurrence 
property, which allows to prove the presence of weakly recurrent points for any 
measurable dynamical system\footnote{In distinction to standard recurrent points 
   which might be absent.}. Since the analysis of the recurrence for the semigroups 
is based on the study of the corresponding property for specially constructed Markov 
chains, we discuss the questions related to the recurrence of measurable Markov chains 
in Section~\ref{s:S2}. Here we follow mainly our recent paper~\cite{Bl2019}, however 
a number of new results (especially about the recurrence of individual trajectories) 
are discussed as well. 
In Section~\ref{s:S3} we give necessary definitions related to the recurrence in 
semigroups and obtain general conditions under which the set of recurrent points 
is the set of full measure for a given reference measure (e.g. the Lebesgue measure). 
Finally Section~\ref{s:S4} is dedicated to the questions of inheritance of the 
recurrence properties from generators to the semigroup and back. In particular, we 
demonstrate that it is possible that the generators might have no recurrent points, 
while almost all points are recurrent for the semigroup. In the backward direction 
it is shown that each generator might have sets of recurrent points of positive 
(Lebesgue) measure, while the semigroup does not have a single uniformly recurrent point.

\section{One generator semigroup -- $S_1$ action.} \label{s:S1} 
We start with the simplest case when the semigroup has a single generator -- 
a measurable mapping $T$ of a Borel compact metric space $(X,\cB)$ into itself. 
Here and in the sequel we always assume that $X$ is equipped with a topology, 
compatible with the $\sigma$-algebra $\cB)$,
and a (non-necessarily dynamically invariant) probabilistic measure $m$, to which 
we will refer as a {\em reference} measure. Let us give 
several definitions related to the notion of recurrence. 

\smallskip

\bdef{We say that a point $x\in X$ is 
\begin{itemize} 
\item {\em recurrent} if for any open neighborhood 
$O_x\ni x$ there exists $t=t(x,O_x)\in\IZ_+$ such that 
$T^tx \in O_x$ (i.e. a trajectory of the point $x$ returns eventually to $O_x$). 
\item {\em weakly recurrent} if for any open neighborhood 
$O_x\ni x$ there exists a point $y\in O_x$ and $t=t(y,O_x)\in\IZ_+$ 
such that $T^ty \in O_x$ (i.e. a trajectory of the point $y$ returns eventually to $O_x$, 
which means that the point $x$ is non-wandering\footnote{A point $x\in X$ is non-wandering 
if for each neighborhood $U$ of $x$, there exists $n>0$ such that $U\cap T^nU\ne\emptyset$.
We introduce the notion of the weak recurrence instead of the non-wandering property 
because its generalization for Markov chains and semigroups is more transparent.}). 
\item {\em uniformly recurrent} if for any open neighborhood 
$O_x\ni x$ the trajectory starting at $x$ returns to $O_x$ infinitely many times. 
\item {\em uniformly weakly recurrent} if for any open neighborhood 
$O_x\ni x$ there exists a point $y\in O_x$ whose trajectory returns to $O_x$ infinitely many times. 
\end{itemize}
}

The set of recurrent points we denote by $\Rec(T)$, while the set of weakly recurrent 
points by $\Rec_w(T)$, and their uniform counterparts by $\Rec_u(T)$ and $\Rec_{wu}(T)$ 
respectively. 

While the second variant describes a weakened version of the recurrence property, 
the last two deal with stronger versions related to the points returning to their neighborhoods  
infinitely many times.

Our main results about properties of various recurrent sets are collected in the following theorem. 

\begin{theorem}\label{t:rec1} 
\begin{itemize}
\item[(0)] $\Rec(T) = \Rec_u(T) \subseteq \Rec_{wu}(T) \subseteq \Rec_w(T) 
                                \supseteq {\rm Clos}(\Rec(T))$.%
\item[(a)] $\forall B\in\cB$ $m$-a.a. points $x\in B$ return to $B$ under 
  dynamics if and only if %
\beq{e:poincare}{  \sum_{n\ge1} m(T^{-n}A) = \infty \quad \forall A\in\cB, ~~m(A)>0} %
\item[(b)] $m(\Rec(T))=1$ iff the property (\ref{e:poincare}) holds true for each 
    {\bf open} $A\in\cB, ~m(A)>0$.
\item[(c)] $\Rec_{wu}(T) \ne \emptyset$.
\end{itemize}
\end{theorem}

Before proceeding to the proof of this result let us make a few comments. 

In the context of the Poincare Theorem it is easy to construct an example when a point returns 
under dynamics to a given set only a finite number of times and wanders away after this. 
In contrast, we show that the simple recurrence of a point implies the uniform recurrence 
(item 0).\footnote{Thanks for the anonymous referee for the unexpected question about this.} 
Unfortunately this implication is wrong for general Markov chains or semigroups.\footnote{There is 
    no control over probabilities of returns to small neighborhoods of a recurrent point.} 
Surprisingly the situation with the weak recurrence property is more 
tricky and the corresponding implication is not valid even for the simple recurrence 
(see example \ref{ex:wu} bellow). 

In general it is possible that $\Rec(T)=\emptyset$, but $\Rec_w(T)$ cannot be empty.  
This follows from to the property (c) above and the trivial inclusion 
$\Rec_w(T) \supseteq \Rec_{wu}(T)$.
On the other hand, there might be weakly recurrent points not belonging to the closure 
of the set of recurrent points (for example, when $\Rec(T)=\emptyset$). 
To demonstrate the last claim consider an example.

\begin{example}\label{ex:1} $X:=[0,1],~\cB:={\rm Bor}$ and 
$$  Tx:=\function{4/5 &\mbox{if~}~~x=0 \\ 
                   x/2 &\mbox{if~}~~0<x\le1/2 \\  
              1-T(1-x) &\mbox{otherwise}.} $$
\end{example}
Under the action of this map any trajectory converges to one of the end-points $0$ or $1$, 
but the trajectory of $0$ converges to $1$ while the trajectory of $1$ converges to $0$. 
Therefore there are no recurrent points but both end-points (and only they) are weakly recurrent.  

From the Poincare Theorem it follows that a necessary condition for the absence of 
recurrent points is the absence of invariant measures. The example~\ref{ex:1} 
is based on the simplest model for this phenomenon, namely on the map %
\beq{e:map}{Tx:=\function{1 &\mbox{if~}~~x=0 \\ 
                   x/2 &\mbox{if~}~~0<x\le1}}
from the unit interval into itself. 
Since in this case $\Rec(T)=\{0\}$, one needs to ``isolate'' the trajectory 
of the origin from the attracting set of this point, which was realized in the example~\ref{ex:1}.

The difference between the formulations of the items (a) and (b) in the theorem above is that 
while in (\ref{e:poincare}) we take into account all measurable sets, only open ones are 
considered in the item (b). To demonstrate the difference\footnote{Thanks for the anonymous 
   referee for this observation.}
observe that for the map~(\ref{e:map}) and the reference measure $m:=\delta_0$ 
we have $T^{-n}\{0\}=\emptyset$ for $n>0$ and hence 
$$ \sum_{n\ge1} m(T^{-n}\{0\}) = 0 < \infty $$ 
for the set $\{0\}$ of full $m$-measure. Thus the property (\ref{e:poincare}) fails 
(in fact here we have a positive measure wandering set  $\{0\}$), 
however the point $0$ is uniformly recurrent and $m(\{0\})=1$. 
 
An important property of the map (\ref{e:map}) is that it  is piecewise contractive. 
Let us analyze up to which extent the latter property is necessary for the absence of 
invariant measures. 

\begin{example}\label{ex:1exp} Let $X$ be the unit circle identified with the 
semi-interval $(0,1]$. Consider a partition of $X$ by the infinite sequence of 
points $\{2^{-n}\}_{n=0}^\infty$ and let $X_k:=(2^{-k},2^{-k+1}],~k\ge1$. 
Then $X$ is a disjoint union of $X_k$. Given $a\in(3/4,3/2)$ define a piecewise 
linear map $T_ax:=ax - (2a-3/2)2^{-k}$ if $x\in X_k$. 
\end{example}

Direct calculation shows that for the given region of values of the parameter $a$ we have 
$$ T_aX_k\in X, \quad T_aX_k\setminus X_k \subset \cup_{i>k}X_i, \quad 
     m(T_aX_k\cap X_k)=m(X_k)/2 ,$$
where $m$ is the Lebesgue measure on $X$. 

In other words the map $T_a$ expands the length of each interval $X_k$ in 
$a$ times and shifts the resulting interval to the left (i.e. to intervals with 
larger indices). Thus any probabilistic measure under the action of $T_a$ 
converges weakly to the Dirac measure at point 1, while 
$T_a\delta_1=\delta_{3/4}\ne\delta_1$. 

An important observation is that the family of maps $T_a$ includes both piecewise 
contractive ($a<1$) and expanding ($a>1$) maps, as well as piecewise 
isometric ($a=1$) maps. 

The point $1$ is recurrent for all maps from this family, but using the map 
in the example~\ref{ex:1exp} as a building block in the same manner as in the 
example~\ref{ex:1}, one gets the map without recurrent points.

\begin{example}\label{ex:wu} Let $X:=[0,1]$ be the unit interval. Consider a partition of $X$ 
by the infinite sequence of points $\{2^{-n}\}_{n=0}^\infty$ and let $X_k:=(2^{-k-1},2^{-k}],~k\ge0$. 
Then $X\setminus\{0\}$ is a disjoint union of $X_k$. Let $T:X\to X$ be defined as follows: %
\beq{e:map-wu}{Tx:=\function{2/3 &\mbox{if~}~~x=0 \\ 
                                               x+2^{-k-2} &\mbox{if~}~~x\in X_k,~ k\ge1 \\
                                               x    &\mbox{otherwise .} }}
\end{example}

In other words the map $T$ shifts each set $X_k$ for $k\ge1$ to the right by the half of 
its length. By the construction $\forall k\ge 1$ he have $TX_k\cap X_k\ne\emptyset$ and 
$T^{n+1}X_k\cap X_{k-i}=\emptyset~~\forall n,i\ge1$. Thus each small enough 
neighborhood $U$ of the point $0$ satisfies the property that $\exists y\in U$ such that 
$Ty\in U$, but each point from $U$ leaves $U$ after a finite number of iterations and 
never returns back. Therefore the point $0$ is weakly recurrent, but not uniformly weakly recurrent. 

\bigskip

\n{\bf Proof} of Theorem~\ref{t:rec1}. The claim about inclusions of various recurrent 
sets is a direct consequence of their definitions, except for the equivalence of the simple 
versions to the uniform ones. To prove this fact consider a nested sequence of neighborhoods 
$O_x^1 \supset O_x^2 \supset \dots \supset O_x^n \supset \dots \ni x$ of a given 
recurrent point $x$. By definition the trajectory of the point $x$ visit each of these 
neighborhoods and thus due to their nested construction does this an infinite number 
of times. In the case of the weak recurrence this argument does not work since 
one needs to consider a trajectory of the point $y$ (which depends on the neighborhood) 
instead of $x$. 

\smallskip

The item (a) is a generalization of the Poincare Theorem 
for an arbitrary reference measure instead of the invariant measure. The proof of this 
result together with the fact that $m(\Rec(T))=1$ if the property (\ref{e:poincare}) holds 
true were obtained recently in our paper \cite{Bl2019} 
and we formulate it here basically for the reader's convenience.

Let us prove item (b). Suppose that (\ref{e:poincare}) holds for each open set of positive 
$m$-measure. Choose a countable basis $\{U^n\}_{n\in\IZ_+}$ of open sets 
in $X$ and set 
$$\t{U}^n:=\{x\in U^n:~T^nx\notin U^n~~\forall n\in\IZ_+\} .$$ 
Then the complement to the union of these sets $\t{U}:=\cup_{n\in\IZ_+}\t{U}^n$ 
is the set of recurrent points. On the other hand, by item (a) $m(\t{U}^n)=0$ for each 
$n\in\IZ_+$. Hence $m(\t{U})=0$, which implies that its complement $\Rec(T)$ is the 
set of full measure.\footnote{The argument used here is similar to the one described 
   in\\ http://planetmath.org/proofofpoincarerecurrencetheorem2}

In the inverse direction we proceed as follows. Choose an arbitrary open set $A\in\cB$ 
with $m(A)>0$. By the assumption $m$-a.a. points from $A$ are recurrent and hence 
are returning to arbitrary small their neighborhoods. On the other hand, since the set 
$A$ is open it contains each point together with some neighborhood. Thus  $m$-a.a. 
points from $A$ return to $A$ under dynamics and hence the inequality \ref{e:poincare} 
holds true.

\smallskip

To prove item (c) we check first that the set of weakly recurrent points is non empty 
for any measurable map, i.e. $\Rec_w(T)\ne\emptyset$.
Assume that there are no weakly recurrent points. For each $x$ we consider 
the function 
$$ R(x):=\sup\{r:~B_r(x) \cap (\cup_{n>0}T^nB_r(x))=\emptyset\} ,$$ 
where $B_r(x):=\{y\in X: \rho(x,y)\le r\}$ -- the ball of radius $r$ centered at the 
point $x$. In words $R(x)$ is the radius of the largest ball around the point $x$ whose 
images do not intersect with this ball. 
Clearly under our assumption on the absence of weakly recurrent points, $R(x)>0$ 
everywhere. Our aim now is to show that $\gamma:=\inf_xR(x)>0$. 
Indeed, let $\{x_n\}$ be an infinite sequence of points in $X$. By the compactness 
of $X$ all limit points of this sequence also belong to $X$ and thus the values 
of the function $R$ at these points are strictly positive.  
Therefore around all points $\{T^nx\}_{n\ge0}$ of the infinite trajectory of the point 
$x\in X$ there are disjoint balls of radius $\gamma$, which contradicts to the 
boundedness of $X$. 

Now we are ready to deal with the similar claim for the smaller set $\Rec_{wu}(T)$. 
Assume from the contrary that $Rec_{wu}(T))=\emptyset$. This means that 
$\forall x\in X~\exists O_x\ni x$ such that $\forall y\in O_x~\exists n_y<\infty$ 
for which $T^ny\notin O_x~\forall n>n_y$. 

Choose a smaller open neighborhood $O'_x$ of a point $x$ such that
$$ {\rm Clos}(O'_x) \subset O_x .$$
For any sequence of points $\{y_k\}\subset O'_x$ each its limit point belongs 
to Clos$(O'_x)\subset O_x$. On the other hand, by the assumption the points from 
$O_x$ return to $O_x$ only finitely many times, which together with the 
compactness of Clos$(O'_x)$ implies that $\sup_k n_{y_k}<\infty$.

Therefore $\forall x\in X~\exists O_x\ni x$ and $n_x<\infty$ such that 
$T^nO_x\cap O_x=\emptyset~\forall n>n_x$. The collection $\{O_x\}_{x\in X}$ 
is an open cover of the compact set $X$. Therefore there exists a finite 
sub-cover $\{C_1,\dots, C_M\}$ for which  
$\forall i~\exists n_i~~T^nC_i\cap C_i=\emptyset$. Let $N:=\max_i n_i<\infty$. 
Hence $T^{kN}C_i\cap C_i=\emptyset~~\forall i, \forall k\in\IZ_+$, which 
cannot happen. The last statement will be proven in a separate combinatorial lemma.

\begin{lemma}\label{l:L1} 
Let $G:\Omega\to2^\Omega\setminus\emptyset=:\t\Omega$ be a multivalued map and 
let $\#\Omega:= M<\infty$. 
Then $\exists \omega\in\Omega$ and $n\in\IZ_+$ such that
$$ G^n\omega \cap \omega \ne \emptyset .$$
\end{lemma}

Since $\Omega\in \t\Omega$, the map $G$ can be extended as a map from the set 
$\t\Omega$ into itself. Therefore all its iterates are well defined.

\proof Choose an arbitrary element $\omega\in\Omega$ and consider 
its trajectory\footnote{A sequence of sets from $\t\Omega$ 
being consecutive images of the one-point set $\{\omega\}$.} 
$\{G^k\omega\}_{k=0}^n$ up to time $n\in\IZ_+$. 
Since the total number of elements in $X$ is finite and each set in the 
trajectory contains at least one element from $\Omega$, we see that whence 
$n>M$ some images of $\omega$ along the trajectory needs to intersect. 
Thus there is a pair of integers $0\le i<j\le M+1$ such that 
$G^i\omega \cup G^j\omega \ne \emptyset$. Then for each element 
$\omega'\in\Omega$ belonging to this intersection we have
$\omega'\in G^{j-i}\omega'$, which proves the claim. \qed

In a sense this result is a simple generalization of the fact 
that each trajectory of a dynamical system with a finite phase space is 
eventually periodic.

\smallskip

This finishes the proof of Theorem~\ref{t:rec1}.\qed

%%%%%%====

\section{Markov chains} \label{s:S2} 
As we already mentioned in the Introduction, a direct generalization of the techniques 
used for the analysis of the recurrence for dynamical system (considered in 
Section~\ref{s:S1}) is not available for semigroups. The main technical tool  
that we use to establish recurrence type results for semigroups is a construction 
of a special Markov chain, whose recurrence properties allow to prove the 
corresponding results for semigroups. Therefore we review here some definitions 
and basic results about Markov chains obtained in our recent paper \cite{Bl2019}. 
To be consistent with the notation of semigroups, we change slightly terminology 
used in that publication. 
In addition, we will obtain a number of new results related to the existence of 
various types of recurrent behavior in Markov chains, which allow us to derive 
the corresponding results for semigroups.

\bdef{By a (homogeneous discrete time) {\em Markov chain} one means a random 
process $\xi_t: (\Omega,\cF,P) \to (X,\cB, m),~t\in\IZ_+$ acting on a Borel 
space $(X,\cB)$ with a finite reference measure $m$ (which needs not to
coincide with the distribution of the process $\xi_t$). The Markov chain (process) $\xi_t$ 
is completely defined by a family of {\em one-step transition probabilities}
$$Q(x,A):=P(\xi_{t+1}\in A|\xi_t=x),~t\in\IZ_+, ~A\in\cB.$$
}

Iterating $Q(\cdot,\cdot)$ one gets a sequence of transition probabilities in $s\in\IZ_+$ time steps:
$$Q^s(x,A):=P(\xi_{t+s}\in A|\xi_t=x),~t\in\IZ_+, ~A\in\cB.$$

Recall also the action of the Markov chains on measures:

$$ Q\mu(A) := \int Q(x,A) d\mu(x), ~~ \forall A\in\cB .$$

\bdef{By the $t$-{\it preimage} with $t\in\IZ_+$ of a set $B\in\cB$
under the action of the Markov chain $\xi_t$ we mean the set of points
$$ Q^{-t}(B) := \{x\in X:~~ Q^t(x,B)>0\} .$$
}
In other words this is the set of initial points of trajectories which 
reach the set $B$ at time $t$ with positive probability.  

Now we are ready to return to the notion of recurrence. 

\bdef{We say that a point $x\in X$ is 
\begin{itemize}
\item {\em recurrent} if for any open neighborhood $O_x\ni x$ there exists 
$t=t(x,O_x)\in\IZ_+$ such that $Q^t(x,O_x)>0$ 
(i.e. a trajectory returns eventually to $O_x$ with positive probability). 
\item {\em weakly recurrent} if for any open neighborhood $O_x\ni x$ 
there exists a point $y\in O_x$ and $t=t(y,O_x)\in\IZ_+$ such that 
$Q^t(y,O_x)>0$ (i.e. a trajectory starting from 
the point $y$ returns eventually to $O_x$ with positive probability). 
\item {\em uniformly recurrent} if 
$\liminf\limits_{n\to\infty} Q^n(x,O_x) > 0$ for any open neighborhood $O_x\ni x$ 
and {\em uniformly weakly recurrent} if 
$\exists y\in O_x~~\liminf\limits_{n\to\infty} Q^n(y,O_x) > 0$. 
\end{itemize}
}
Similarly to the deterministic setting by $\Rec(Q), \Rec_w(Q), \Rec_u(Q), \Rec_{wu}(Q)$ 
we denote the sets of recurrent, weakly recurrent, uniformly recurrent and weakly 
uniformly recurrent points respectively.

Comparing these definitions with their deterministic counterparts
(see e.g. \cite{Ni,Si}) or to the notions of recurrence and
transience well studied for the case of countable Markov chains
(see e.g. \cite{Fe,MT,Or,Sp}) one is tempted to make the
conditions stronger assuming that the corresponding events take
place with probability one (instead of just being positive).
Examples discussed in \cite{Bl2019} demonstrate absence of the recurrence 
under this stronger assumption even in the simplest situations. 
As we will see, similar effects take place in the case of semigroups as well.
Reasons why in the definitions of various types of recurrence we assume $Q^n(x,A)>0$ 
instead of $=1$ are discussed in detail in \cite{Bl2019}.\footnote{The main among them 
is that typically under the stronger assumption the set of recurrent points shrinks dramatically.} 
In Example~\ref{ex:Q-u} we will show that in general $\Rec(Q) \ne \Rec_u(Q)$.
On the other hand, denoting by $\Rec^s(Q), \Rec_u^s(Q)$ the sets of (uniformly) recurrent 
points for which the eventual return occurs with probability 1, one can show (using exactly the 
same argument as in the proof of Theorem~\ref{t:rec1} (item 0) ) that $\Rec^s(Q)=\Rec_u^s(Q)$.

Since a dynamical system $(T,X)$ generates a Markov chain by means of the transition 
probabilities $Q(x,A):=1_{A}(Tx)$, the example \ref{ex:1} demonstrates that a general 
Markov chain needs not possess even a single recurrent point. 

Our main results about properties of various recurrent sets in the probabilistic 
setting may be formulated as follows.

\begin{theorem}\label{t:Markov-full} 
\begin{itemize}
\item[(0)] $\Rec(Q) \supseteq \Rec_u(Q) \subseteq \Rec_{wu}(Q) \subseteq \Rec_w(Q) 
                             \supseteq {\rm Clos}(\Rec(Q))$.
\item[(a)]  If %
\beq{e:rec-Q}{ \sum_{n\ge1} m(Q^{-n}(A) \cap A) =\infty~~\forall A\in\cB: m(A)>0 ,}%
then $m(\Rec(Q))=1$ and $\Rec_w(Q)\supseteq{\rm supp}(m)$.
\item[(b)] Rec$_w(Q)\ne\emptyset$.
\item[(c)] If $\mu$ is a probabilistic invariant measure for the Markov chain 
defined by the transition probabilities $Q(\cdot,\cdot)$, then $\mu(\Rec(Q))=1$. 
\end{itemize}
\end{theorem}

\proof The claim about inclusions of various recurrent sets is a direct consequence 
of their definitions.\footnote{See also example~\ref{ex:Q-u}, which demonstrates the 
    non-equivalence of the sets $\Rec(Q)$ and $\Rec_u(Q)$. Compare this to a comment 
    about the equality $\Rec^s(Q)=\Rec_u^s(Q)$ above.  }

The part about the recurrence in item (a) was proven in Theorems 4,5 of \cite{Bl2019}, 
while the part about the weak recurrence is a consequence of the inclusions between 
the recurrence sets.  

Technically the proof of this result is based on the generalization of the classical 
Poincare Recurrence theorem, which says that if the measure $m$ is dynamically 
invariant, then for any measurable set $m$-almost all points from this set return 
eventually under dynamics. A generalization of this claim for Markov chains 
instead of dynamical systems is not very difficult. Another point is how to deal 
with non-dynamically invariant measures $m$. See \cite{Bl2019} for details. 

The item (b) is a probabilistic version of Theorem~\ref{t:rec1}(c). 
Unfortunately the argument used in the proof of Theorem~\ref{t:rec1}(c) cannot 
be used in this more complicated setting and we develop a different approach, 
also based on the compactness of the phase space.

Assume that there are no weakly recurrent points. Therefore for each point 
$x\in X$ there exists an open neighborhood $O_x$ such that 
$$ Q^n(y,O_x)=0~~  \forall y\in O_x, \forall n\in\IZ_+ .$$ 
The collection of neighborhoods $\{O_x\}_{x\in X}$ is a cover by open sets of the 
compact $X$. Hence there exists a finite subcover $\{C_1,\dots,C_M\}$, which by 
the construction should satisfy the property:
$$ Q^n(x,C_i)=0~~\forall i\in\{1,\dots,M\}, x\in C_i, n\in\IZ_+ .$$ 
It remains to show that this property cannot hold for $n$ large enough. 

Consider the set $\Omega:=\cup_{i=1}^M C_i$ and define a map $G$ as follows: 
$$ G(C_i) := \cup_{j=1}^M\{C_j:~~\exists x\in C_i: ~~Q(x,C_j)>0 \}.$$ 
Applying Lemma~\ref{l:L1} to this setting we observe that there exists $C_i$ from our 
open cover and a positive integer $n$ such that $G^n(C_i)\cap C_i\ne\emptyset$, 
which contradicts to the construction of the subcover $\{C_i\}$ and hence to the 
original assumption about the absence of weakly recurrent points. 

In item[c] we consider the Markov chain with the stationary distribution $\mu$ 
and thus one can introduce the corresponding Markov shift (see, e.g. \cite{Ve}). 
Thus we get a measurable map (Markov shift) with the invariant measure for 
which all results from the previous section can be applied. In particular, from 
Theorem~\ref{t:rec1}(b) we get the desired statement. \qed

At the moment we do not have conditions under which the sets of (weakly) 
uniformly recurrent points for a Markov chain are sets of full measure for 
a general reference measure $m$. 
However using the argument already applied in the proof of item (c) one gets 
this result in the case when the reference measure 
is invariant with respect to the Markov chain.

In the sequel we will need to compare recurrence properties of different Markov chains. 

\bdef{We say that transition probabilities $Q(\cdot,\cdot)$ and $\t{Q}(\cdot,\cdot)$ 
are {\em compatible} if for each $x\in X$ the measures $Q(x,\cdot)$ and $\t{Q}(x,\cdot)$ 
are absolutely continuous with respect to each other; and {\em uniformly compatible} if 
\beq{e:compat}{
   0 < \gamma\t{Q}(x,A) \le Q(x,A) \le \frac1\gamma \t{Q}(x,A)~~\forall x\in X, A\in\cB .}
}

Clearly the uniform compatibility implies the compatibility.

\begin{theorem}\label{t:c4} 
Let $Q(\cdot,\cdot)$ and $\t{Q}(\cdot,\cdot)$ be (uniformly) compatible transition 
probabilities of two Markov chains defined on the same space. 
Then $\Rec(Q)=\Rec(\t{Q})$ and $\Rec_w(Q)=\Rec_w(\t{Q})$. 
\end{theorem}
\proof We prove a slightly more general property, namely that for a given pair of a 
point $x\in X$ and a measurable set $A\subset X$ the probability to reach $A$ starting 
from the point $x$ under the actions of the Markov chains $Q$ and $\t{Q}$ is either 
positive in both cases or zero in both cases. 

Let us show that the compatibility of $Q(\cdot,\cdot)$ and $\t{Q}(\cdot,\cdot)$ 
implies the compatibility of $Q^n(\cdot,\cdot)$ and $\t{Q}^n(\cdot,\cdot)$ for all $n\ge2$. 
Indeed, by definition of the transition probabilities  
$$ Q^{n+1}(x,A) = \int Q^n(y,A)Q(x,dy) .$$
This gives the absolute continuity of the the measures $Q^2(x,\cdot)$ and $\t{Q}^2(x,\cdot)$ 
with respect to each other (it is enough to consider the corresponding integral sums). 
In the case of uniform compatibility we get
$$ 0 < \gamma^2\t{Q}^2(x,A) \le Q^2(x,A) \le \gamma^{-2} \t{Q}^2(x,A)~~\forall x\in X, A\in\cB .$$
Continuing this calculations, by the standard induction argument we get the desired statement. 
 
Therefore for each $n\in\IZ_+, x\in X, A\in\cB$ either $Q^n(x,A) + \t{Q}^n(x,A) =0$ 
or $Q^n(x,A) \cdot \t{Q}^n(x,A) >0$.
Thus, choosing an arbitrary neighborhood $O_x$ of the point $x$ and setting $A:=O_x$ 
we get the claim. \qed

\begin{remark} It is worth mentioning that the measures $Q^n(x,\cdot)$ and $Q^k(x,\cdot)$ 
for $n\ne k$ need not to be absolutely continuous with respect to each other (and typically are not). 
\end{remark}

\begin{remark} Unfortunately in order to apply Theorem~\ref{t:c4} to check the recurrence 
of a certain point $x$ under the action of two different Markov chains one needs to assume 
that both chains are compatible everywhere, rather not locally near the point $x$. 
\end{remark}

\section{General finitely generated semigroups -- $S_d$ action} \label{s:S3} 

\bdef{By a finitely generated {\em semigroup} $\cS$ of maps one means 
a collection of endomorphisms from a set $X$ into itself admitting 
a finite sub-collection $\bT:=\{T_1,\dots,T_d\}\subset \cS$ satisfying the 
condition that each element from $\cS$ may be represented in the form %
\beq{e-repr}{ \cS \ni T_{[i_1,\dots,i_n]} := T_{i_n}\circ T_{i_{n-1}}\circ\dots\circ T_{i_1},
    \qquad i_k \in\{1,\dots,d\},\quad k\in\{1,\dots,n\}  .}
The sub-collection $\bT$ is called the semigroup set of {\em generators}. 
If there are no generators $T_i$ that cannot be represented in the form 
(\ref{e-repr}) without $T_i$, then the generators set is said to be {\em minimal}.
}

In what follows we will not make any additional assumptions about the 
correspondence between the generators, which means that we 
consider only {\em free} semigroups. 

Naturally a finitely generated semigroup is completely characterized by the choice of 
its generators, however there might be several sets of generators, and  one needs to 
check that this choice does not change properties of the semigroup under study. 
Further on we will show that all our results do not depend on a particular representation 
of a semigroup by the set of generators $\bT$.  

To simplify the notation for a given set of generators $\bT$, we also refer to the 
action of the semigroup $\cS$ as the action of $\bT$. 

\bdef{By a {\em trajectory} of length $n$ composed by the applications of the 
generators $\bT$ of the semigroup $\cS$ starting from a point $x\in X$ one means 
an ordered collection of points 
$$\{T_{i_1}(x), T_{[i_1,i_2]}(x),\dots,  T_{[i_1,\dots,i_n]}(x)\} ,$$
were the indices $i_k\in\{1,2,\dots,d\}$. }

Denote by $N(x,A,n)$ the number of trajectories of length $n$ of the semigroup $\bT$, 
starting from the point $x\in X$ and ending in the set $A\in\cB$. 

\bdef{We say that a point $x\in X$ is 
\begin{itemize}
\item {\em recurrent} if for each open neighborhood $O_x$ of a point $x\in X$ 
we have $N(x,O_x,n) > 0$ for some $n=n(x)<\infty$.
\item {\em weakly recurrent} if for each open neighborhood 
$O_x$ of a point $x\in X$ there is a point $y\in O_x$ such that
$N(y,O_x,n) > 0$ for some $n=n(y)<\infty$.
\item {\em uniformly recurrent} if for each open neighborhood 
$O_x$ of a point $x\in X$ the proportion of trajectories finishing at the set $O_x$ is positive, i.e. 
$$ \kappa(x,O_x):= \liminf_{n\to\infty} \frac{N(x,O_x,n)}{N(x,X,n)} > 0 .$$
\item {\em weakly uniformly recurrent} if for each open neighborhood 
$O_x$ of a point $x\in X$ there is a point $y\in O_x$ such that
$$ \kappa(y,O_x):= \liminf_{n\to\infty} \frac{N(y,O_x,n)}{N(y,X,n)} > 0 .$$
\end{itemize}
}

In the last two versions the strong condition requires positivity of a ``share'' of returning 
trajectories starting from the point $x$, while the weak condition corresponds to the 
existence of a point nearby with a positive ``share'' of returning trajectories.

In the same manner as in the previous sections by $\Rec(\bT)$, $\Rec_w(\bT)$, $\Rec_u(\bT)$ 
and $\Rec_{wu}(\bT)$ we denote the sets of recurrent, weakly recurrent, uniformly recurrent 
and weakly uniformly recurrent points respectively.

Similarly to the case of Markov chain instead of the assumption $\kappa(x,O_x)>0$ one 
can use a stronger version $\kappa(x,O_x)=1$. Denote by $\Rec^s(\bT), \Rec_u^s(\bT)$ 
the sets of strong (uniformly) recurrent points under the latter assumption.\footnote{This 
   in general leads to a dramatical shrinkage of the corresponding recurrent sets.} 
In Example~\ref{ex:Q-u} we will show that in general $\Rec(\bT) \ne \Rec_u(\bT)$.
On the other hand, using connections between the semigroups and Markov chains and the 
discussion immediately before Teorem~\ref{t:Markov-full} one shows that 
$\Rec^s(\bT) = \Rec_u^s(\bT)$.

In order to study the recurrence property for semigroups it is natural to try to generalize 
the condition (\ref{e:poincare}) for this more general setting. Unfortunately a naive direct 
generalization of type %
\beq{e:poincare2}{\sum_{i=1}^d  \sum_{n=0}^\infty m(T_i^{-n}A) = \infty 
                            ~~{\rm if}~~ m(A)>0 }
obviously does not work. As a counterexample assume that one of the maps transforms the entire 
phase space to a set of zero $m$-measure (see Section~\ref{s:S4} for details).
A more promising generalization is as follows. Choose a sequence of indices 
$\{i_k\}_{k=1}^\infty$ with $i_i\in\{1,\dots,d\}$. Then one might expect that 
a non-autonomous version of (\ref{e:poincare}) might work, namely that %
\beq{e:poincare3}{\sum_{n=1}^\infty m(T_{[i_1,i_n]}^{-1}A) = \infty 
                            ~~{\rm if}~~ m(A)>0 }
implies that under the action of $T_{[i_1,\dots]}$ $m$-almost all points from the  
set $A$ return back to $A$. Nevertheless a close look to this assumption shows 
that even this is not correct. The point is that the inheritance of the recurrence property 
from individual generators to the semigroup and back is a rather delicate problem. 
We will discuss this in detail in Section~\ref{s:u-rec}, where several counterexamples 
demonstrating quite counterintuitive recurrence properties will be constructed and studied. 

To overcome this difficulty we propose a construction based on the presentation 
of a semigroup as a Markov chain acting on the same phase space.

For a given set of generators $\bT$ and a probability distribution $\bp:=\{p_1,\dots,p_d\}$ 
consider a Markov chain defined by the following transition probabilities: %
\beq{e:Q}{ Q(x,A):= \sum_i^d p_i1_A(T_ix) .}%

Clearly, for a given minimal set of generators $\bT$ and a non-degenerate distribution 
$\bp:=\{p_1,\dots,p_d\}$ (i.e. $\prod_ip_i>0$) there is a bijection between the set 
of trajectories of the semigroup generated by $\bT$ and the set of realizations 
of the Markov chain $Q$. 

\begin{theorem}\label{t:T=Q} Let $\bT$ be a generator set, $\bp:=\{p_1,\dots,p_d\}$ 
be a non-degenerate distribution, and let $Q(\cdot,\cdot)$ be the family of transition 
probabilities for the Markov chain defined by the relation (\ref{e:Q}). 
Then 
\begin{itemize}
\item [(a)] $\Rec(\bT)=\Rec(Q)$ and  $\Rec_w(\bT)=\Rec_w(Q)$.
\item [(b)] Let additionally $p_1=p_2=\dots=p_d$ and let $m$ be the invariant measure 
for the Markov chain $Q$, then $m(\Rec_u(\bT))=m(\Rec_{wu}(\bT))=1$.
\end{itemize}
\end{theorem}

\begin{corollary}
\begin{itemize}
\item[(i)] $\Rec(\bT) \supseteq \Rec_u(\bT) \subseteq \Rec_{wu}(\bT) 
                      \subseteq \Rec_w(\bT) \supseteq {\rm Clos}(\Rec(\bT))$.
\item[(ii)] If $Q$ satisfies the condition~(\ref{e:rec-Q}) for some non-degenerate distribution 
$\bp$, then $m(\Rec(\bT))=1$ and $\Rec(\bT)\supseteq\supp(m)$. 
\end{itemize}
\end{corollary}

\n{\bf Proof} of Theorem~\ref{t:T=Q}. 
Let $\bp$ be a uniform distribution over the set of indices $\{1,\dots,d\}$. 
Then for any $x\in X$, each its neighborhood $O_x$ and any positive integer $n$ we have 
$$ \kappa(x,O_x,n):= \frac{N(x,O_x,n)}{N(x,X,n)} = Q^n(x,O_x) .$$
Therefore for this specific choice of the distribution $\bp$ whenever $x\in\Rec(Q)$ we have 
that $\kappa(x,O_x,n)>0$ for each neighborhood $O_x$ and some $n<\infty$. 
Hence $N(x,O_x,n)>0$ and thus  $x\in\Rec(\bT)$. 
Similarly we have the inverse inclusion, which proves the claim in item~(a) about the set of 
recurrent points for the uniform distribution $\bp$. 

The part about the weakly recurrent points follows from the same argument. 

\smallskip

To deal with a general non-degenerate but non-necessarily uniform distribution $\bp$ 
one needs to compare properties of Markov chains generated by the same semigroup's 
generators but with different distributions $\bp$.   

\begin{lemma}\label{l:probs} 
Let the distributions $\bp:=\{p_1,\dots,p_d\}$ and $\t\bp:=\{\t{p}_1,\dots,\t{p}_d\}$ 
be non-degenerate. Then the corresponding Markov chains $Q$ and $\t{Q}$ satisfy 
the inequalities (\ref{e:compat}) and thus are uniformly compatible. 
\end{lemma}
\proof This claim follows from the simple observation that
$$ 0<\min_i\frac{p_i}{\t{p}_i }\le \frac{\sum_i^d p_i1_A(T_ix) }{\sum_i^d \t{p}_i1_A(T_ix) } 
                                         \le \max_i\frac{p_i}{\t{p}_i} < \infty $$
due to the non-degeneracy of the distributions $\bp$ and $\t\bp$. \qed

Applying now the result of Theorem~\ref{t:c4} for the Markov chains with 
uniformly compatible transition probabilities we finish the proof of the claim 
in item (a) of Theorem~\ref{t:T=Q}. 

\smallskip

To prove the claim in item (b) we apply the construction of the Markov chain 
with the transition probabilities (\ref{e:Q}) (exactly as in the proof of item (a)) 
and  the uniform distribution $\bp$. Observing that under this assumption 
$\kappa(x,O_x,n)=Q^n(x,O_x)$ we see that all information about the statistics 
of values $\kappa(x,O_x,n)$ can be derived from the statistics of $Q^n(x,O_x)$. 
Therefore to get the result one applies the claim of Theorem~\ref{t:Markov-full}(c). 

This completes the proof of Theorem~\ref{t:T=Q}. \qed

To demonstrate that in general $\Rec(\bT)\ne\Rec_u(\bT)$ consider the following example. 

\begin{example}\label{ex:Q-u} Let $X:=[0,1]$ be the unit interval and let $\bT$ be 
generated by a pair of maps $T_i:X\to X$ defined by the relations $T_1x:=x^2, ~T_2x:=1$.
\end{example}

In this example $\Rec(\bT)=\{0,1\} = \Rec_w(\bT) \ne \Rec_u(\bT)=\Rec_{wu}(\bT)=\{1\}$. 

Consider now the Markov chain, generated by the semigroup $\bT$ in this example with 
equal probabilities $p_1=p_2=1/2$, i.e. by (\ref{e:Q}) the corresponding transition 
probabilities can be written as $Q(x,A):=\frac12 1_{A}(x^2) + \frac12 1_{A}(1)$. 
Thus we get an example of a Markov chain for which the simple recurrence does not 
imply the uniform one.

\bigskip

As we already mentioned the same semigroup may be generated  by different 
generator sets. Therefore we need to check that the sets of various recurrent 
points of a semigroup do not depend on its representation by different 
generator sets.

\begin{theorem}\label{t:representation} Let $\bT:=\{T_1,\dots,T_d\}$ and 
 $\t\bT:=\{\t{T}_1,\dots,\t{T}_{\t{d}}\}$ be two minimal generator sets of the 
same semigroup of maps. Then $\Rec(\bT)=\Rec(\t\bT)$ and $\Rec_w(\bT)=\Rec_w(\t\bT)$.
\end{theorem}
\proof By the definition of the minimal generator set $\forall j$ we have the following 
presentation %
\beq{e:pres}{ \t{T}_j = T_{[i_1,i_2,\dots,i_n]} ,} %
where $i_k\in\{1,\dots,d\}$. The number of elements in the minimal generator 
set is finite and hence there exists $N<\infty$ such that each $\t{T}_j\in\t\bT$ is composed 
by at most $N$ generators from $\bT$.

Choose an arbitrary non-degenerated distribution $\bp:=\{p_1,\dots,p_d\}$. Then accordingly 
to (\ref{e:pres}) we construct a new distribution $\t\bp:=\{\t{p}_1,\dots,\t{p}_{\t{d}}\}$ 
by the relations
$$ \t{p}_j:=\prod_{k=1}^{n_j} p_{j_k} ,$$
where $n_j\le N~~\forall j$. Hence the distribution $\t\bp$ is again non-degenerated. 

Therefore it remains to check that the Markov chains with transition probabilities 
$Q$ and $\t{Q}$ are compatible. To this end we apply the same argument as 
in the proof of Lemma~\ref{l:probs} with a very slight modification related to 
the different (but well controlled by means of the constant $N$ above) numbers 
of addends in the nominator and the denominator. \qed

\section{Counterexamples and recurrence inheritance} \label{s:S4} 

This section is dedicated to the problem of inheritance of the recurrence property 
from the generators of the semigroup to the semigroup itself and back. 

We assume again that we have a metric space $X$ equipped with a Borel $\sigma$-algebra 
$\cB$ of measurable subsets and on this space we consider the action of a semigroup 
with a finite set of generators $\bT:=\{T_1,\dots,T_d\}$. 

\subsection{Simple recurrence}\label{s:s-rec}

The next result demonstrates that the direct inheritance of the simple recurrence property 
is rather straightforward.

\begin{lemma}\label{t:inc} $\cup_i \Rec(T_i) \subseteq \Rec(\bT)$.
\end{lemma}
\proof Assume that $x \in \Rec(T_i)$. By the definition of the recurrent point for 
a single map we have that for each neighborhood $O_x\ni x ~ \exists 0<n_x<\infty$ 
such that $T_i^{n_x}x\in O_x$. 
Hence $N(x, O_x, n_x) \ge 1 >0$, which proves that $\Rec(T_i) \subseteq \Rec(\bT)$. \qed

Surprisingly the inverse inheritance is more tricky. 
Assume that $\Rec(\bT)\ne\emptyset$. The following example shows that 
despite this assumption it is possible that none of the generators 
of the semigroup $\bT$ has recurrent points. Moreover, the set of recurrent points 
of the semigroup might be everywhere dense while each semigroup generator 
possess no recurrent points.

\begin{example}\label{ex:2} Let $X:=[0,1]$ be the unit interval. 
Consider a pair of maps $T_i:X\to X$ defined by the relations  %
\beq{e:diag1}{
    T_1x:=\function{\frac12x+\frac14 &\mbox{if } 0\le x < \frac12 \\ 
                             \frac12x+\frac12 &\mbox{if } \frac12\le x < 1\\
                              0 & \mbox{if } x=1,} \quad
    T_2x:=\function{\frac12x &\mbox{if } 0< x \le \frac12 \\ 
                             \frac12x+\frac14 &\mbox{if } \frac12 < x \le 1\\
                              1 & \mbox{if } x=0}
}
\end{example}

\begin{figure}
\begin{center}\begin{tikzpicture}
\put(-80,0){\draw  (0,0) -- (4,0) -- (4,4) -- (0,4) -- cycle; \draw (2,0) to (2,4); \draw (0,2) to (4,2);
                      \draw [->,line width=1.2pt] (0, 1) to (2,2); \draw [->,line width=1.2pt] (2, 3) to (4,4);
                      \draw [fill] (4,0) circle (0.08cm); \node at (0.5,3.5){$T_1$};
                    }
\put(80,0){\draw  (0,0) -- (4,0) -- (4,4) -- (0,4) -- cycle; \draw (2,0) to (2,4); \draw (0,2) to (4,2);
                      \draw [<-,line width=1.2pt] (0, 0) to (2,1); \draw [<-,line width=1.2pt] (2, 2) to (4,3);
                      \draw [fill] (0,4) circle (0.08cm); \node at (0.5,3.5){$T_2$};
                    }
 \end{tikzpicture} \end{center}
\caption{Graphs of the maps $T_i$ in the Example~\ref{ex:2}}\label{f:zigzag} \end{figure}
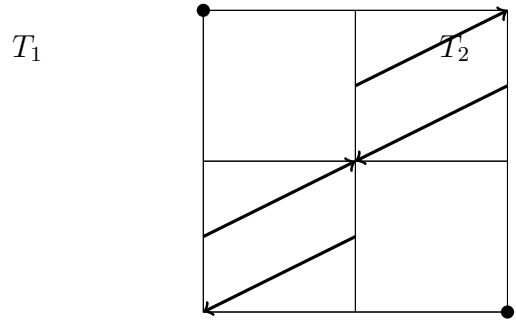

The graphs of the maps $T_i$ are shown on the Figure~\ref{f:zigzag}. 
Under the action of $T_1$ all the points from the interval $[0,1/2)$ converge 
to the point $1/2$, but $T_11/2=3/4$, and all the points from the interval $[1/2,1)$ 
converge to the point $1$, while $T_11=0$. 
Similarly under the action of $T_2$ all the points from the interval $(0,1/2]$ converge 
to the point $0$, but $T_20=1$, and all the points from the interval $(1/2,1]$ 
converge to the point $1/2$, while $T_21/2=1$. 
Therefore neither of these maps have recurrent points.

To study $\Rec(\bT)$ for $\bT:=\{T_1,T_2\}$ we consider the corresponding Markov 
chain with the uniform probability distribution $p_1=p_2=1/2$. We claim that this 
Markov chain has exactly two ergodic Lebesgue invariant measures supported by the 
intervals $[0,1/2]$ and $[1/2,1]$ correspondingly, which implies that the closure of 
$\Rec(\bT)$ coincides with $[0,1]$. To check this claim consider an even simpler 
Markov chain $\t{Q}$ on $X:=[0,1]$, defined by the transition probabilities 
$\t{Q}(x,A):=\frac12(1_A(\frac{x}2+\frac12) + 1_A(\frac{x}2))$ and corresponded  
to the simplest linear iteration function system. Such systems are well studied and 
it is known (see,  e.g. \cite{Bl2001}) that $\t{Q}$ has the only one ergodic invariant 
measure which is equal to the Lebesgue measure on $[0,1]$. 

Observe now that the construction of the Example~\ref{ex:2} consists of two building 
blocks, where the (shifted and normalized) maps $\frac{x}2+\frac12$ and $\frac{x}2$ 
are applied. Therefore the Markov chain corresponding to this example has 
exactly two ergodic Lebesgue invariant measures supported by the intervals $[0,1/2]$ 
and $[1/2,1]$ correspondingly. Thus the closure of $\Rec(\bT)$ coincides with $[0,1]$. 

\subsection{Uniform recurrence}\label{s:u-rec}

The inheritance of the recurrence property becomes much more involved when 
we are interested in statistics of trajectories of recurrent points in the semigroups 
under study. The claim of Lemma~\ref{t:inc} is no longer valid for the uniform 
recurrence setting. In particular, the next example demonstrates that $\Rec_u(\bT)$ might 
be empty even if each of generators possesses uniformly recurrent points.

\begin{example}\label{ex:3} Let $X:=[0,1]$ be the unit interval. 
Consider a pair of maps $T_i:X\to X$ defined by the relations  %
\beq{e:diag2}{
    T_1x:=\function{\frac12x &\mbox{if } 0\le x < \frac13 \\ 
                             \frac14 & \mbox{if } x=0,\\
                             \frac12x+\frac13 &\mbox{if } \frac13< x < \frac12\\
                             \frac35 & \mbox{if } x=\frac13,\\
                             1-T_1(1-x) & \mbox{otherwise} }
                            \qquad
    T_2x:=\function{x + \frac13 &\mbox{if } 0\le x < \frac13 \\ 
                            -\frac12x+\frac34 &\mbox{if } \frac13 \le x < \frac12\\
                              \frac13 & \mbox{if } x=\frac12\\
                             1-T_1(1-x) & \mbox{otherwise} }
}
\end{example}

\begin{figure}
\begin{center}\begin{tikzpicture}
\put(-90,0){\draw  (0,0) -- (6,0) -- (6,6) -- (0,6) -- cycle; 
                  \draw (2,0) to (2,6); \draw (4,0) to (4,6);   \draw (0,2) to (6,2); \draw (0,4) to (6,4);
                      \draw [<-,line width=1.2pt] (0, 0) to (2,1); \draw [->,line width=1.2pt] (4, 5) to (6,6);
                      \draw [fill] (0,1.6) circle (0.08cm); \draw [fill] (6,4.4) circle (0.08cm); 
                      \draw [<->,line width=1.2pt] (2, 2) to (3,2.5); 
                      \draw [<-,line width=1.2pt] (4, 4) to (3,3.5); 
                      \draw [fill] (2,3.6) circle (0.08cm); \draw [fill] (4,2.4) circle (0.08cm); 
                      \node at (0.5,3.5){$T_1$};
                    }
\put(90,0){\draw  (0,0) -- (6,0) -- (6,6) -- (0,6) -- cycle; 
                  \draw (2,0) to (2,6); \draw (4,0) to (4,6);   \draw (0,2) to (6,2); \draw (0,4) to (6,4);
                      \draw [->,line width=1.2pt] (0, 2) to (2,4); \draw [<-,line width=1.2pt] (4,2) to (6, 4);
                      %\draw [fill] (0,1.6) circle (0.08cm); \draw [fill] (6,4.4) circle (0.08cm); 
                      \draw [<-,line width=1.2pt] (3, 3) to (2,3.5); 
                      \draw [<-,line width=1.2pt] (3, 3) to (4,2.5); 
                      \draw [fill] (3,2) circle (0.08cm);
                      \node at (0.5,3.5){$T_2$};
                    }
 \end{tikzpicture} \end{center}
\caption{Graphs of the maps $T_i$ in the Example~\ref{ex:3}}\label{f:no-u-rec} \end{figure}
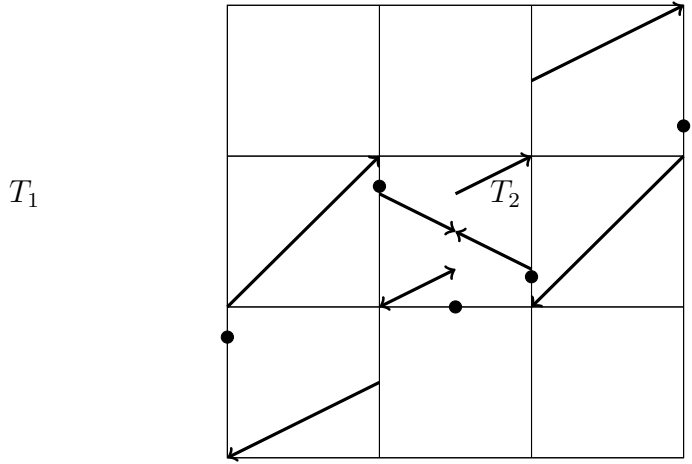

The graphs of the maps $T_i$ are shown on the Figure~\ref{f:no-u-rec}. 

The map $T_1$ transforms each of the intervals $(0,1/3), (1/3,2/3), (2/3,1)$ into itself 
and the points $0,1/3,2/3,1$ are uniformly recurrent. 
The map $T_2$ transforms each of the intervals $(0,1/3)$ and $(2/3,1)$ onto the interval 
$(1/3,2/3)$, while the latter interval is mapped into itself. For this map $\Rec_u(T_2)=\{1/2\}$. 

To study uniformly recurrent points of the semigroup generated by these maps, consider 
the corresponding Markov chain with equal probabilities $p_1=p_2=1/2$. 
Observe that the only candidates for being uniformly recurrent points for the Markov chain are 
the points $0,1/2,1$. On the other hand, a trajectory may converge to one of these points 
only if a long enough series of consecutive applications of the same map $T_1$ or $T_2$ 
is applied. The probability of a realization of such an event vanishes exponentially with the 
length of the series. Therefore $Q^n(x,O_x)$ goes to zero as $n\to\infty$ for each small 
enough neighborhood $O_x$. Thus $\Rec_u(\bT)=\emptyset$.

\bigskip

It is natural to ask if this unexpected inheritance is due to the fact that each of the generators 
has only a few recurrent points? To this end, without any changes in the arguments above 
one can assume that the map $T_1$ is identical on the intervals $[0,1/3)$ and $(2/3,1]$, which 
gives $m(\Rec(T_1))=2/3$. On the other hand, one cannot do anything of this sort with the 
second generator $T_2$. In order to overcome this restriction we consider yet another example.

\begin{example}\label{ex:4} Let $X:=[0,1]$ be the unit interval. 
Consider a pair of maps $T_i:X\to X$ defined by the relations  %
\beq{e:diag}{
    T_1x:=\function{x &\mbox{if } 0\le x < \frac13 \\ 
                             \frac12x+\frac13 &\mbox{if } \frac13< x < \frac12\\
                             \frac12x+\frac29 & \mbox{if } \frac12< x < \frac23\\
                             x-\frac13 & \mbox{if } \frac23< x \le 1 }
                            \qquad
    T_2x:=\function{x + \frac13 &\mbox{if } 0\le x < \frac13 \\ 
                            -\frac12x+\frac34 &\mbox{if } \frac13 \le x \le \frac23\\
                              x & \mbox{if } \frac23< x \le 1  }
}
\end{example}

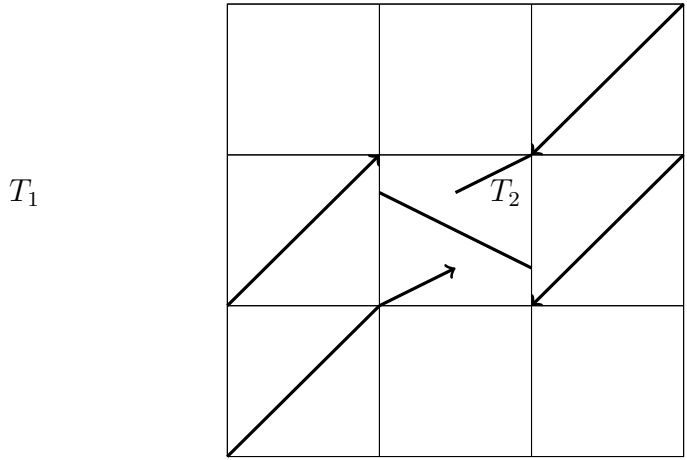
\begin{figure}
\begin{center}\begin{tikzpicture}
\put(-90,0){\draw  (0,0) -- (6,0) -- (6,6) -- (0,6) -- cycle; 
                  \draw (2,0) to (2,6); \draw (4,0) to (4,6);   \draw (0,2) to (6,2); \draw (0,4) to (6,4);
                      \draw [line width=1.2pt] (0, 0) to (2,2);
                      \draw [->,line width=1.2pt] (2, 2) to (3,2.5); 
                      \draw [line width=1.2pt] (4, 4) to (3,3.5); 
                      \draw [->,line width=1.2pt] (6, 4) to (4,2); 
                      \node at (0.5,3.5){$T_1$};
                    }
\put(90,0){\draw  (0,0) -- (6,0) -- (6,6) -- (0,6) -- cycle; 
                  \draw (2,0) to (2,6); \draw (4,0) to (4,6);   \draw (0,2) to (6,2); \draw (0,4) to (6,4);
                      \draw [->,line width=1.2pt] (0, 2) to (2,4); \draw [<-,line width=1.2pt] (4,4) to (6, 6);
                      \draw [line width=1.2pt] (2, 3.5) to (4,2.5);
                      \node at (0.5,3.5){$T_2$};
                    }
 \end{tikzpicture} \end{center}
\caption{Graphs of the maps $T_i$ in the Example~\ref{ex:4}}\label{f:no-u-rec-Leb} \end{figure}

The graphs of the maps $T_i$ are shown on the Figure~\ref{f:no-u-rec-Leb}. 

In this example $\Rec(T_1)=[0,1/3]$ and $\Rec(T_2)=\{1/2\}\cup[2/3,1]$, so 
$m(\Rec(T_1))=m(\Rec(T_2))=1/3$. On the other hand, by the same arguments as in the 
analysis of the example~\ref{ex:3}, there are no uniformly recurrent points for the semigroup 
generated by the maps $T_1,T_2$. 

\bigskip

To demonstrate that in general the inverse inheritance also does not take place 
it is enough to consider the Example~\ref{ex:2}.

\end{document}